\documentclass[psamsfonts,english]{amsart} 

    \usepackage[T1]{fontenc}
    \usepackage{babel}                          
    \usepackage{natbib}                         
    \usepackage{amsmath}                        
    \usepackage{amssymb}
    \usepackage{amsthm}                         
    \usepackage[mathcal]{eucal}                 
    \usepackage{url}

    \newcommand{\lilleskrift}{\fontsize{8pt}{1}}

    \newtheorem{teo}{Theorem}
    \theoremstyle{definition}
    \newtheorem{defn}[teo]{Definition}
    \newtheorem{ex}{Example}
    \theoremstyle{remark}
    \newtheorem{rem}[teo]{Remark}

    \newcommand{\FF}{\mathbb{F}}
    \newcommand{\QQ}{\mathbb{Q}}
    \newcommand{\ZZ}{\mathbb{Z}}
    \newcommand{\RR}{\mathbb{R}}
    
    \newcommand{\PP}{\mathbb{P}}
    \newcommand{\jac}[1]{\mathcal{J}_{#1}}
    \newcommand{\heltal}[1]{\mathfrak{O}_{#1}}

    \DeclareMathOperator{\End}{End}
    \DeclareMathOperator{\Div}{Div}
    \DeclareMathOperator{\Mat}{Mat}
    \DeclareMathOperator{\gal}{Gal}

\begin{document}

\title[Large Cyclic Subgroup of Jacobians of HEC]
{Large Cyclic Subgroups of Jacobians of Hyperelliptic Curves}


\author[C.R. Ravnshøj]{Christian Robenhagen Ravnshøj}

\address{Department of Mathematical Sciences \\
Faculty of Science \\
University of Aarhus \\
Ny Munkegade \\
Building 1530 \\
DK-8000 Aarhus C}

\email{cr@imf.au.dk}

\thanks{Research supported in part by a Ph.D. grant from CRYPTOMAThIC}

\keywords{Jacobians, hyperelliptic curves, complex multiplication, cryptography}

\subjclass[2000]{Primary 14H40; Secondary 11G15, 14Q05, 94A60}


\begin{abstract}
In this paper we obtain conditions on the divisors of the group order of the Jacobian of a hyperelliptic genus~$2$ curve,
generated by the complex multiplication method described by \cite{weng03} and \cite{gaudry}. Examples,
where these conditions imply that the Jacobian has a large cyclic subgroup, are given.
\end{abstract}

\maketitle

\section{Introduction}

In elliptic curve cryptography it is essential to know the number of points on the curve. Cryptographically we are
interested in curves with large cyclic subgroups. Such elliptic curves can be constructed. The construction is based on the
theory of complex multiplication, studied in detail by \cite{atkin-morain}. It is referred to as the \emph{CM method}.

\cite{koblitz89} suggested the use of hyperelliptic curves to provide larger group orders. Therefore constructions of
hyperelliptic curves are interesting. The CM method for elliptic curves has been generalized to hyperelliptic curves of genus~$2$ by
\cite{spallek}, and efficient algorithms have been proposed by \cite{weng03} and \cite{gaudry}.

Both algorithms take as input a primitive, quartic CM field $K$, and give as output a hyperelliptic genus~$2$ curve $C$ over a
prime field $\FF_p$. A prime number $p$ is chosen such that $p=\omega\overline\omega$ for a number
$\omega\in\heltal{K}$, where $\heltal{K}$ is the ring of integers of $K$. We have $K=\QQ(\eta)$ and
$K\cap\RR=\QQ(\sqrt{D})$, where $\eta=i\sqrt{a+b\xi}$ and
    $$\xi=\begin{cases}
    \frac{1+\sqrt{D}}{2}, & \textrm{if $D\equiv 1\pmod{4}$,} \\
    \sqrt{D}, & \textrm{if $D\equiv 2,3\pmod{4}$}.
    \end{cases}
    $$
Write $\omega=c_1+c_2\xi+(c_3+c_4\xi)\eta$, $c_i\in\ZZ$. Let $C$ be a hyperelliptic curve of genus~$2$ over $\FF_p$ with
$\End(C)\simeq\heltal{K}$. The Jacobian $\jac{C}(\FF_p)$ is isomorphic to
    \begin{equation}\label{eq:jac(Fp)(0)}
    \ZZ/n_1\ZZ\times\ZZ/n_2\ZZ\times\ZZ/n_3\ZZ\times\ZZ/n_4\ZZ,
    \end{equation}
where $n_i\mid n_{i+1}$ and $n_2\mid p-1$. In this paper, conditions on the prime divisors of the number $n_2$ are
obtained, and examples, where these conditions imply that the Jacobian $\jac{C}(\FF_p)$ has a large cyclic
subgroup, are given. The conditions on the prime divisors are given by the following theorem.

{\samepage
\begin{teo}\label{teo:maal}
Let $C/\FF_p$ be a hyperelliptic curve of genus~$2$ with $\End(C)\simeq\heltal{K}$, where $K$ is a primitive, quartic CM field. Assume that the structure of $\jac{C}(\FF_p)$ is given by \eqref{eq:jac(Fp)(0)}. Let $\ell\mid n_2$ be an odd prime number. Then $\ell\leq Q$, where
    \begin{align*}
    Q &= \max\{a,D,a^2-b^2D\}, \\
    \intertext{if $D\equiv 2,3\pmod{4}$, and}
    Q &= \max\{a,D,4a(a+b)-b^2(D-1),aD+2b(D-1)\},
    \end{align*}
if $D\equiv 1\pmod{4}$. If $\ell>D$, then $c_1\equiv 1\pmod{\ell}$ and $c_2\equiv 0\pmod{\ell}$.
\end{teo}
}

\begin{rem}
Since the number $n_2\mid p-1$ and $\ell\mid n_2$, it follows that $\ell\neq p$.
\end{rem}

\section{Hyperelliptic curves}

A hyperelliptic curve is a smooth, projective curve $C\subseteq\PP^n$ of genus $g\geq 2$ with a separable, degree $2$
morphism $\phi:C\to\PP^1$. Let $C$ be a hyperelliptic curve of genus $g=2$ defined over a prime field $\FF_p$, where
$\FF_p$ is of characteristic~$p>2$. By the Riemann-Roch theorem there exist an embedding $\psi:C\to\PP^2$, mapping $C$
to a curve given by an equation of the form
    $$y^2=f(x),$$
where $f\in\FF_p[x]$ is of degree $\deg(f)=6$ and have no multiple roots \cite[see][chapter~1]{cassels}.

The set of principal divisors $\mathcal{P}(C)$ on $C$ constitutes a subgroup of the degree 0 divisors $\Div_0(C)$. The
Jacobian $\jac{C}$ of $C$ is defined as the quotient
    $$\jac{C}=\Div_0(C)/\mathcal{P}(C).$$
Let $\ell\neq p$ be a prime number. The $\ell^n$-torsion subgroup $\jac{C}[\ell^n]<\jac{C}$ of elements of order
dividing $\ell^n$ is then by \cite[theorem~6, p.~109]{lang59}
    $$\jac{C}[\ell^n]\simeq\ZZ/\ell^n\ZZ\times\ZZ/\ell^n\ZZ\times\ZZ/\ell^n\ZZ\times\ZZ/\ell^n\ZZ.$$
An endomorphism $\varphi:\jac{C}\to\jac{C}$ induces a $\ZZ_\ell$-linear map
    $$\varphi_\ell:T_\ell(\jac{C})\to T_\ell(\jac{C})$$
on the $\ell$-adic Tate-module $T_\ell(\jac{C})$ of $\jac{C}$ \cite[chapter~VII, \S1]{lang59}. Hence
$\varphi$ is represented on $\jac{C}[\ell]$ by a matrix $M\in\Mat_{4\times 4}(\ZZ/\ell\ZZ)$. Let $P(X)\in\ZZ[X]$ be the characteristic polynomial of $\varphi$ \cite[see][pp.~109--110]{lang59} and $P_M(X)\in(\ZZ/\ell\ZZ)[X]$ the characteristic polynomial of $M$. Then \cite[theorem~3, p.~186]{lang59}
    \begin{equation}\label{eq:KarPolKongruens}
    P(X)\equiv P_M(X)\pmod{\ell}.
    \end{equation}
Since $C$ is defined over $\FF_p$, the mapping $(x,y)\mapsto (x^p,y^p)$ is an isogeny on~$C$. This isogeny induces an
endomorphism $\varphi$ on the Jacobian $\jac{C}$, the Frobenius endo\-morphism. The characteristic polynomial $P(X)$ of
$\varphi$ is of degree $4$ \cite[theorem~2, p.~140]{tate}. Theorem~\ref{teo:maal} will be established by using the
identity~\eqref{eq:KarPolKongruens} on the Frobenius.


\section{CM fields}\label{sec:CMfields}

An elliptic curve $E$ with $\ZZ\neq\End(E)$ is said to have \emph{CM}. Let $K$ be an ima\-ginary, quadratic number field
with ring of integers $\heltal{K}$. $K$ is a \emph{CM field}. If \mbox{$\End(E)\simeq\heltal{K}$}, then $E$ is said to
have \emph{CM by $\heltal{K}$}. More generally a CM field is defined as follows.

\begin{defn}[CM field]
A number field $K$ is a CM field, if $K$ is a totally imaginary, quadratic extension of a totally real number field
$K_0$.
\end{defn}

In this paper only CM fields of degree $[K:\QQ]=4$ are considered. Such a field is called a \emph{quartic} CM field. Let
$K_0=K\cap\RR$. Then $K_0$ is a real, quadratic number field, $K_0=\QQ(\sqrt{D})$. Since $K$ is a totally imaginary,
quadratic extension of $K_0$, a number $\eta\in K$ exists, such that $K=K_0(\eta)$, $\eta^2\in K_0$. The number $\eta$
is totally imaginary,
and we may assume $\eta=i\eta_0$, $\eta_0\in\RR$, and that $-\eta^2$ is totally positive.

Let $C$ be a hyperelliptic curve of genus $g=2$. Then $C$ is said to have CM by~$\heltal{K}$, if
$\End(C)\simeq\heltal{K}$. The structure of $K$ determines whether $C$ is irreducible. More precisely, the following
theorem holds.

\begin{teo}\label{teo:reducibel}
Let $C$ be a hyperelliptic curve of genus~$2$ with CM by $\heltal{K}$, where $K$ is a quartic CM field. Then $C$ is
reducible if, and only if, $K/\QQ$ is Galois with Galois group $\gal(K/\QQ)\simeq\ZZ/2\ZZ\times\ZZ/2\ZZ$.
\end{teo}

\begin{proof}
\cite[proposition~26, p.~61]{shi}.
\end{proof}

Theorem~\ref{teo:reducibel} motivates the following definition.

\begin{defn}[Primitive, quartic CM field]
A quartic CM field $K$ is called primitive if either $K/\QQ$ is not Galois, or $K/\QQ$ is Galois with cyclic Galois
group.
\end{defn}

\section{The CM method for genus~$2$}

The CM method for genus~$2$ is described in detail by \cite{weng03} and \cite{gaudry}. In short, the CM method is based
on the construction of the class polynomials of the number field $K$. The prime number $p$ has to be chosen such that
    \begin{equation}\label{eq:norm}
    p=\omega\overline\omega
    \end{equation}
for a number $\omega\in\heltal{K}$. There are 2 approaches to choose such a prime number~$p$. Either pick a random
prime number $p$, and try to solve the complex norm equation~\eqref{eq:norm} in $\heltal{K}$, or generate a number
$\omega\in\heltal{K}$, such that $\omega\overline\omega$ is a prime number. The first approach needs deep theory, e.g. class groups. The second can be implemented in a short algorithm, and is based on elementary theory. Moreover,
empirical results indicate that the elementary method is the faster of the two approaches \cite[tab\-le~1]{weng03}. Thus the elementary method
is preferable. The algorithm is given in figure~\ref{fig:p-konstruktion} for $D\equiv 2,3\pmod{4}$. The algorithm for
$D\equiv 1\pmod{4}$ is similar \cite[section~8]{weng03}.

\begin{rem}\label{rem:primiske}
In either way we get an $\omega\in\heltal{K}$ with $\omega\overline\omega=p$. We may assume that $\omega$ fulfils the
additional condition $\gcd(c_3,c_4)=1$, where the numbers $c_3$ and~$c_4$ are given by equation~\eqref{eq:omega} in
section~\ref{sec:properties}. In the first approach, if $\omega$ does not fulfil this condition, we can just pick
another prime number $p$. In the elementary method we can incorporate this condition in the algorithm.
\end{rem}

\begin{figure}[bt]
\begin{description}
    \item[Input] CM-field $K=\QQ\left(i\sqrt{a+b\sqrt{D}}\right)$.
    \item[Output] Prime $p=\omega\overline\omega$ and $\omega\in\heltal{K}$.
    \begin{enumerate}
    \item Choose random numbers $c_3,c_4\in\ZZ$ such that $\gcd(c_3,c_4)=1$ and $c_3^2b-c_4^2bD\equiv 0\pmod{2}$.
    \item Set $2n:=-2c_3c_4a-c_3^2b-c_4^2bD$.
    \item Choose $c_1$ at random as a divisor of $n$.
    \item Set $c_2:=n/c_1$.
    \item Set $p:=c_1^2+c_2^2D+c_3^2a+c_4^2aD+2c_3c_4bD$. If $p$ is not a prime number, start again.
    \item Set $\omega:=c_1+c_2\sqrt{D}+(c_3+c_4\sqrt{D})i\sqrt{a+b\sqrt{D}}$.
    \end{enumerate}
\end{description}
\caption{Elementary method to choose a prime number $p=\omega\overline\omega$ in the case $D\equiv
2,3\pmod{4}$.}\label{fig:p-konstruktion}
\end{figure}

\section{Properties of $\jac{C}(\FF_p)$}\label{sec:properties}

Let $K$ be a primitive, quartic CM field with real subfield $K_0=\QQ(\sqrt{D})$ of class number $h(K_0)=1$. Write
$K=\QQ(\eta)$, where $\eta=i\sqrt{a+b\xi}$ and
    $$\xi=\begin{cases}
    \frac{1+\sqrt{D}}{2}, & \textrm{if $D\equiv 1\pmod{4}$,} \\
    \sqrt{D}, & \textrm{if $D\equiv 2,3\pmod{4}$}.
    \end{cases}
    $$
We may assume that $a\pm b\sqrt{D},a+b\frac{1\pm\sqrt{D}}{2}>0$, cf. section~\ref{sec:CMfields}. Let $p$ be a prime
number such that
    $$p=\omega\overline\omega$$
for a number $\omega\in\mathfrak{O}=\heltal{K_0}+\eta\heltal{K_0}$. Since $h(K_0)=1$, we can write
    \begin{equation}\label{eq:omega}
    \omega=c_1+c_2\xi+(c_3+c_4\xi)\eta,\quad c_i\in\ZZ.
    \end{equation}
We may assume $\gcd(c_3,c_4)=1$, cf. remark~\ref{rem:primiske}. Let $C/\FF_p$ be a hyperelliptic curve of genus~$2$ with
CM by $\heltal{K}$. Write
    \begin{equation}\label{eq:jac(Fp)}
    \jac{C}(\FF_p)\simeq\ZZ/n_1\ZZ\times\ZZ/n_2\ZZ\times\ZZ/n_3\ZZ\times\ZZ/n_4\ZZ,
    \end{equation}
where $n_i\mid n_{i+1}$ and $n_2\mid p-1$ \cite[see][proposition~5.78, p.~111]{hhec}. Depending on the remainder of $D$ modulo $4$, we obtain conditions on the prime divisors of the number
$n_2$.

{\samepage
\begin{teo}\label{teo}
Let $C/\FF_p$ be a hyperelliptic curve of genus~$2$ with CM by $\heltal{K}$. Assume that the structure of $\jac{C}(\FF_p)$ is
given by~\eqref{eq:jac(Fp)}. Let $\ell\mid n_2$ be an odd prime number. Then $\ell\leq Q$, where
    \begin{align*}
    Q &= \max\{a,D,a^2-b^2D\}, \\
    \intertext{if $D\equiv 2,3\pmod{4}$, and}
    Q &= \max\{a,D,4a(a+b)-b^2(D-1),aD+2b(D-1)\},
    \end{align*}
if $D\equiv 1\pmod{4}$. If $\ell>D$, then $c_1\equiv 1\pmod{\ell}$ and $c_2\equiv 0\pmod{\ell}$.
\end{teo}
}

\begin{proof}
Assume $D\equiv 2,3\pmod{4}$. Since $\omega\overline\omega=p$ we find that
    \begin{align}
    p &= c_1^2+c_2^2D+c_3^2a+c_4^2aD+2c_3c_4bD,  \label{eq:p} \\
    0 &= 2c_1c_2+c_3^2b+c_4^2bD+2c_3c_4a.        \label{eq:0}
    \end{align}
Let $P(X)$ be the characteristic polynomial of the Frobenius $\varphi$.
    $$
    P(X) = \prod_{i=1}^4(X-\omega_i)
         = X^4-4c_1X^3+(2p+4(c_1^2-c_2^2D))X^2-4c_1pX+p^2.
    $$
Here $\omega_i$ are the roots of $P(X)$.

Let $\ell\mid n_2$ be an odd prime number. Then by equation~\eqref{eq:jac(Fp)} the Jacobian $\jac{C}(\FF_p)$ contains a
subgroup $U\simeq(\ZZ/\ell\ZZ)^3$. As
    $$(\ZZ/\ell\ZZ)^3<\jac{C}(\FF_p)[\ell]<\jac{C}[\ell],$$
the Frobenius $\varphi$ is represented on $\jac{C}[\ell]$ by a matrix
    $$
    M = \begin{bmatrix}
        1 & 0 & 0 & m_1 \\
        0 & 1 & 0 & m_2 \\
        0 & 0 & 1 & m_3 \\
        0 & 0 & 0 & m_4
        \end{bmatrix}
    $$
Notice that $m_4=\det(M)\equiv\deg(\varphi)=p^2\pmod{\ell}$. Since $p\equiv 1\pmod{\ell}$, $M$ has the characteristic
polynomial
    $$P_M(X) \equiv (X-1)^4 = X^4-4X^3+6X^2-4X+1\pmod{\ell}.$$
Now $P(X)\equiv P_M(X)\pmod{\ell}$. Thus
    $$c_1\equiv c_1^2-c_2^2D\equiv 1\pmod{\ell},$$
since $\ell\neq 2$.

Assume $\ell>D$. Then
    \begin{equation}\label{eq:c1c2}
    c_1\equiv 1\pmod{\ell},\quad c_2\equiv 0\pmod{\ell}.
    \end{equation}
By the equations \eqref{eq:p} and \eqref{eq:0}, we get
    \begin{align*}
    c_1^2+c_2^2D+c_3^2a+c_4^2aD+2c_3c_4bD &\equiv 1 \pmod{\ell}, \\
    2c_1c_2+c_3^2b+c_4^2bD+2c_3c_4a &\equiv 0 \pmod{\ell}.
    \end{align*}
Therefore, by equation~\eqref{eq:c1c2}, the following holds.
    \begin{align}
    c_3^2a+c_4^2aD+2c_3c_4bD &\equiv 0 \pmod{\ell}, \label{eq:c3c4}\\
    c_3^2b+c_4^2bD+2c_3c_4a &\equiv 0 \pmod{\ell}. \nonumber
    \end{align}
It follows that
    $$c_3c_4(a^2-b^2D)\equiv 0\pmod{\ell}.$$
Here $a^2-b^2D=(a+b\sqrt{D})(a-b\sqrt{D})>0$, since $a\pm b\sqrt{D}>0$. Assume $\ell>a^2-b^2D$. Then we get
$c_3c_4\equiv 0\pmod{\ell}$. Thus either $c_3\equiv 0\pmod{\ell}$ or $c_4\equiv 0\pmod{\ell}$.

Assume $\ell>a$. If $c_3\equiv 0\pmod{\ell}$, then $c_4^2aD\equiv 0\pmod{\ell}$ by equation~\eqref{eq:c3c4}, i.e.
$c_4\equiv 0\pmod{\ell}$. On the other hand if $c_4\equiv 0\pmod{\ell}$, then $c_3^2a\equiv 0\pmod{\ell}$, i.e.
$c_3\equiv 0\pmod{\ell}$.

Summing up, $c_3\equiv c_4\equiv 0\pmod{\ell}$ if $\ell>\max\{a,D,a^2-b^2D\}$. But this contradicts $\gcd(c_3,c_4)=1$.
Therefore $\ell\leq \max\{a,D,a^2-b^2D\}$, and the case $D\equiv 2,3\pmod{4}$ is established.

\newpage

Now consider the case $D\equiv 1\pmod{4}$. Since $\omega\overline\omega=p$, we now find that
    \begin{align*}
    p = {} & c_1^2+c_1 c_2+\frac{1}{4}c_2^2(1+D) +c_3^2\Big(a+\frac{1}{2}b\Big)+c_3c_4\Big(\frac{1}{2}b(D+1)+a\Big) \\
           & +c_4^2\Big(\frac{1}{8}b(3D+1)+\frac{1}{4}a(D+1)\Big), \\
    0 = {} & c_1c_2+\frac{1}{2}c_2^2+\frac{1}{2}c_3^2b+c_3c_4(a+b)+c_4^2\Big(\frac{1}{8}b(D+3)+\frac{1}{2}a\Big).
    \end{align*}
The characteristic polynomial of the Frobenius $\varphi$ is given by
    \begin{align*}
    P(X) = {} & X^4-(4c_1+2c_2)X^3+(2p+(2c_1+c_2)^2-c_2^2D)X^2 \\
              & -(4c_1+2c_2)pX+p^2.
    \end{align*}

Let $\ell\mid n_2$ be an odd prime number. As in the case $D\equiv 2,3\pmod{4}$, the Frobenius $\varphi$ is represented on
$\jac{C}[\ell]$ by a matrix $M$ with the characteristic polynomial
    $$P_M(X) \equiv X^4-4X^3+6X^2-4X+1\pmod{\ell}.$$
Since $P(X)\equiv P_M(X)\pmod{\ell}$, it follows that
    $$4c_1+2c_2\equiv (2c_1+c_2)^2-c_2^2D \equiv 4 \pmod{\ell}.$$

Assume $\ell>D$. Then
    $$c_1\equiv 1\pmod{\ell},\quad c_2\equiv 0\pmod{\ell}.$$
Now
    \begin{align*}
    c_3^2(8a+4b)+c_3c_4\left(4b(D+1)+8a\right) & \\
    +c_4^2\left(b(3D+1)+2a(D+1)\right) &\equiv 0 \pmod{\ell} \\
    4c_3^2b+8c_3c_4(a+b)+c_4^2\left(b(D+3)+4a\right) &\equiv 0 \pmod{\ell}.
    \end{align*}
Therefore
    \begin{align}
    4c_3^2a+2c_3c_4b(D-1)+c_4^2(a+b)(D-1) &\equiv 0 \pmod{\ell}, \label{eq:c3} \\
    4c_3^2b+8c_3c_4(a+b)+c_4^2\left(b(D+3)+4a\right) &\equiv 0 \pmod{\ell}. \nonumber
    \end{align}
It follows that
    $$(b^2(D-1)-4a(a+b))(2c_3c_4-c_4^2)\equiv 0\pmod{\ell}.$$
Notice that
    $$4a(a+b)-b^2(D-1)=4\left(a+b\frac{1+\sqrt{D}}{2}\right)\left(a+b\frac{1-\sqrt{D}}{2}\right)>0.$$
Now assume $\ell>4a(a+b)-b^2(D-1)$. Then
    $$2c_3c_4-c_4^2\equiv 0\pmod{\ell}.$$
Thus either $c_4\equiv 0\pmod{\ell}$ or $c_4\equiv 2c_3\pmod{\ell}$.

Assume $\ell>a$. If $c_4\equiv 0\pmod{\ell}$, then $c_3^2\equiv 0 \pmod{\ell}$ by equation~\eqref{eq:c3}, i.e.
$c_3\equiv 0\pmod{\ell}$. This contradicts $\gcd(c_3,c_4)=0$. So $c_4\not\equiv 0\pmod{\ell}$. Then $c_4\equiv
2c_3\pmod{\ell}$. From equation~\eqref{eq:c3} it follows that
    $$c_4^2(2b(D-1)+aD)\equiv 0\pmod{\ell},$$
i.e. $c_4\equiv 0\pmod{\ell}$ if $\ell>2b(D-1)+aD$. But then $c_3\equiv c_4\equiv 0\pmod{\ell}$, a~contradiction.
\end{proof}

\begin{rem}
The condition $\gcd(c_3,c_4)=1$ may be relaxed. In the proof of theorem~\ref{teo}, we only need
$\ell\nmid\gcd(c_3,c_4)$.
\end{rem}

\section{Examples}

By theorem~\ref{teo}, large prime divisors of the order $N=|\jac{C}(\FF_p)|$ will not divide the divisor $n_2$ of $N$.
This is useful if we want to determine the possible cyclic subgroups of $\jac{C}(\FF_p)$.

\begin{ex}
In $K=\QQ\left(i\sqrt{2+\sqrt{2}}\right)$, the prime number
    $$p=15314033922152826237436247359259334919$$
is the complex norm of the number
    \begin{align*}
    \omega = {} & 3913314953099587393-31\sqrt{2} \\
                & +(4483312578+6978049007\sqrt{2})i\sqrt{2+\sqrt{2}}.
    \end{align*}
The CM method yields a hyperelliptic genus~$2$ curve $C$ with Jacobian of order
    \lilleskrift
    $$N=234519634968847474692278544362349582158321382804023011720188699330496198748.$$
    \normalsize
Since $N=2^2\cdot 7^3\cdot 17\cdot 23\cdot 4993\cdot r$, where
    \lilleskrift
    $$r=87556173808919520163329861675989739433243040373597074857097140343$$
    \normalsize
is a prime number, either
    $$\jac{C}(\FF_p)\simeq \ZZ/N\ZZ \quad \textrm{or} \quad \jac{C}(\FF_p)\simeq \ZZ/n_3\ZZ\times\ZZ/n_4\ZZ,$$
where $n_3\in\{2,7,14\}$.
\end{ex}

\begin{ex}
In $K=\QQ\left(i\sqrt{7+\sqrt{5}}\right)$, the prime number
    $$p=14304107096878940330893123933$$
is the complex norm of the number
    \begin{align*}
    \omega = {} & -119599766860084+5279155\sqrt{5} \\
                & +\left(13860963299+4898901569\sqrt{5}\right)i\sqrt{7+\sqrt{5}}.
    \end{align*}
The CM method yields a hyperelliptic genus~$2$ curve $C$ with Jacobian of order
    \lilleskrift
    $$N=204607479838989309536748148297333557447111046976589088984.$$
    \normalsize
Since $N=2^3\cdot 7^3\cdot 71\cdot r$, where
    \lilleskrift
    $$r=1050217015557576630891205130257738047915611254140091$$
    \normalsize
is a prime number, either
    $$\jac{C}(\FF_p) \simeq \ZZ/n_3\ZZ\times\ZZ/n_4\ZZ,$$
where $n_3\in\{1,2,7,14\}$, or
    $$\jac{C}(\FF_p) \simeq \ZZ/2\ZZ\times\ZZ/n_3\ZZ\times\ZZ/n_4\ZZ,$$
where $n_3\in\{2,14\}$.
\end{ex}

\end{document}